\documentclass[twocolumn]{article}
\usepackage[utf8]{inputenc}
\usepackage[left=2cm, right=2cm, top=2cm]{geometry} 
\usepackage[small]{titlesec}
\usepackage{etoolbox}
\usepackage{amsthm}
\makeatletter
\patchcmd{\ttlh@hang}{\parindent\z@}{\parindent\z@\leavevmode}{}{}
\patchcmd{\ttlh@hang}{\noindent}{}{}{}
\makeatother

	\usepackage{mathptmx} % assumes new font selection scheme installed
	\usepackage{times} % assumes new font selection scheme installed
	\usepackage{amsmath} % assumes amsmath package installed
	\usepackage{amssymb}  % assumes amsmath package installed
	\usepackage{algorithm}
	\usepackage{algpseudocode}
	\usepackage{color}
	\usepackage{bm}
	\usepackage{makecell}
	\usepackage{subfig, float}
	\usepackage{tikz,pgfplots}
	\usepackage{epstopdf}
	\usepackage{natbib}
	\usepackage{mathtools}

	\usetikzlibrary{external,arrows,shapes,positioning,plotmarks, arrows.meta}
	\usepgfplotslibrary{external} 
	\newlength\figureheight 
	\newlength\figurewidth  
	
	\DeclareMathAlphabet{\mathcalOld}{OMS}{cmsy}{m}{n}
	\DeclareMathOperator*{\minimizer}{minimize}
	\DeclareMathOperator*{\subjecttoo}{subject\:to}
	\DeclareMathOperator*{\find}{find}

	\newcommand{\alphabm}{\bm{\alpha}}
	\newcommand{\nubm}{\bm{\nu}}
	\newcommand{\etabm}{\bm{\eta}}
	\newcommand{\nbm}{\bm{n}}
	\newcommand{\xbm}{\bm{x}}
	\newcommand{\ubm}{\bm{u}}
	\newcommand{\taubm}{\bm{\tau}}
	\newcommand{\sship}{\mathcalOld{S}_b}
	\newcommand{\vship}{\mathcalOld{V}_b}
	\newcommand{\sspatial}{\mathcalOld{S}_{\mathrm{env}}}
	\newcommand{\vspatial}{\mathcalOld{V}_{\mathrm{env}}}
	\newcommand{\espatial}{\mathcalOld{E}_{\mathrm{env}}}
	\newcommand{\gspatial}{\mathcalOld{G}_{\mathrm{env}}}

	\newcommand{\eworld}{\mathcalOld{E}_{\mathrm{obst}}}

\title{\LARGE \bf
An Optimization-Based Motion Planner for Autonomous Maneuvering of Marine Vessels in Complex Environments \footnote{This work was partially supported by FFI/VINNOVA and the Wallenberg Artificial Intelligence, Autonomous Systems and Software Program (WASP) funded by Knut and Alice Wallenberg Foundation.}
}

\date{}
\author{Kristoffer Bergman \thanks{Kristoffer Bergman, Oskar Ljungqvist and Daniel Axehill are with the Division of Automatic Control,
		Link{\"o}ping University, Link\"oping, Sweden
		(e-mail: {\tt\small \{kristoffer.bergman, oskar.ljungqvist, daniel.axehill\}@liu.se})} \and Oskar Ljungqvist \footnotemark[2] \and Jonas Linder \thanks{Jonas Linder is at ABB AB, Corporate Research, V\"aster{\aa}s, Sweden (e-mail:
		\tt\small jonas.x.linder@se.abb.com)} \and Daniel Axehill \footnotemark[2]  
}% <-this % stops a space 

\begin{document}
	
	% declaration of the new block
	\algblock{ParFor}{EndParFor}
	% customising the new block
	\algnewcommand\algorithmicparfor{\textbf{parfor}}
	\algnewcommand\algorithmicpardo{\textbf{do}}
	\algnewcommand\algorithmicendparfor{\textbf{end\ parfor}}
	\algrenewtext{ParFor}[1]{\algorithmicparfor\ #1\ \algorithmicpardo}
	\algrenewtext{EndParFor}{\algorithmicendparfor}

	\maketitle
	\thispagestyle{empty}
	\pagestyle{empty}

	%%%%%%%%%%%%%%%%%%%%%%%%%%%%%%%%%%%%%%%%%%%%%%%%%%%%%%%%%%%%%%%%%%%%%%%%%%%%%%%%		
		\textbf{\textit{Abstract ---}}\textbf{The task of maneuvering ships in confined environments is a difficult task for a human operator.  
			One major reason is due to the complex and slow dynamics of the ship which need to be accounted for in order to successfully steer the vehicle. 
			In this work, a two-step optimization-based motion planner is proposed for autonomous maneuvering of ships in constrained environments such as harbors. 
			A lattice-based motion planner is used in a first step to compute a feasible, but suboptimal solution to a discretized version of the motion planning problem. 
			This solution is then used to enable efficient warm-start and as a terminal manifold for a second receding-horizon improvement step. 
			Both steps of the algorithm use a high-fidelity model of the ship to plan feasible and energy-efficient trajectories. 
			Moreover, a novel algorithm is proposed for automatic computation of spatial safety envelopes around the trajectory computed by the lattice-based planner. 
			These safety envelopes are used in the second improvement step to obtain collision-avoidance constraints which complexity scales very well with an increased number of surrounding obstacles. 
			The proposed optimization-based motion planner is evaluated with successful results in a simulation study for autonomous docking problems in a model of the Cape Town~harbor.     }
	
\section{Introduction} \label{sec:intro}
Autonomous ship technology is currently witnessing an increased interest in order to improve efficiency and safety of future sea transports. 
Historically, complex and precise maneuvering of ships in harbors, e.g., docking, required external help from support vessels such as tug boats. 
In recent years, as modern ships have been equipped with more thrusters with better capabilities, the ships' ability to independently perform advanced maneuvers has increased. 
Meanwhile, to aid humans operating modern ships equipped with high degrees of redundancy in their truster setup, various control systems have been developed. 
These systems are designed to alleviate the need of manually allocating thruster commands~\cite{fossen2011handbook}, and to automatically maintain a desired course and set speed or vehicle pose (\emph{dynamic positioning})~\cite{perez2006ship,veksler2016dynamic}. 
As a consequence, the current trends in advanced motion control of ships are moving towards an increased degree of automation. 

Although a large amount of different feedback-control techniques has been proposed in the literature for marine vessels, see, e.g.,~\cite{fossen2011handbook,veksler2016dynamic,perez2006ship,moreira2007path}, there only exists a limited amount of work that focus on simplified versions of the trajectory-planning problem in complex and constrained environments, e.g.,~\cite{martinsen2019autonomous,du2018motion}. 
Compared to agile unmanned surface vehicles (USVs), where maritime collision avoidance maneuvers can be computed by partially or completely neglecting the system dynamics~\cite{han2020autonomous,kufoalor2019autonomous}, reliable maneuvering of ships requires a motion planner that accounts for the ship's dynamics and physical constraints in order to compute feasible, safe and optimized trajectories~\cite{martinsen2019autonomous}. 

In this work, a two-step optimization-based motion planning framework is proposed for automatic maneuvering of ships in cluttered environments such as harbors. The framework is based on techniques presented in our recent work~\cite{bergman2019improved,bergman2019bimproved,bergman2020optimization} and is here tailored and extended to account for the specific needs related to ships.   
A lattice-based motion planning algorithm~\cite{pivtoraiko2009differentially} is used in a first step to compute a resolution optimal solution to a discretized version of the motion planning problem. 
By employing classical graph-search algorithms, the lattice-based planner concatenates motion segments from a finite library of optimized maneuvers~\cite{bergman2019improved} and is responsible for solving combinatorial aspects of the motion planning problem, such as selecting which side to pass an obstacle.
Due to its deterministic properties and efficiency, lattice-based motion planners have been used with great success for a variety of different vehicular systems~\cite{pivtoraiko2009differentially,LjungqvistCDC2018,LjungqvistJFR2019,likhachev2009planning,cirillo2014lattice}. Similar motion planning techniques have also been developed for USVs~\cite{du2018motion} using a library of basic curve segments. 
However, since a lattice-based planner is restricted to use a discretized search space, the computed solutions often suffer from discretization artifacts~\cite{andreasson2015fast}. 

To improve the solution computed by the lattice-based planner, a second optimization-based improvement step is used that is warm-started with the trajectory computed by the lattice-based planner, similar to~\cite{bergman2019bimproved}. 
The improvement step can either optimize the complete trajectory at once, or it can run in a receding-horizon fashion over a sliding time-window as proposed in~\cite{bergman2020optimization}. 
In the latter approach, the output computed by the lattice-based planner is used both for enabling an efficient warm-start and as a terminal manifold to ensure feasibility and convergence to the desired goal state. 

The proposed receding-horizon improvement step shares similarities with the work in~\cite{martinsen2019autonomous}, but instead of using a quadratic terminal cost to represent a possibly under-estimate of the remaining cost-to-go after the prediction horizon, the nominal solution computed by the lattice-based planner provides the receding horizon optimization-step with an over-estimate of the remaining cost-to-go. This means that feasibility can be ensured also in complex environments~\cite{bergman2020optimization}. 
Furthermore, compared to~\cite{martinsen2019autonomous} which restricts the ship's body to a manually selected convex polytope, a novel collision-avoidance algorithm is proposed to automatically compute spatial safety envelopes 
along each sample of the nominal solution computed by the lattice-based planner. 
These safety envelopes are described, e.g., using convex polytopes, and represent the constraints imposed on the ship's body from surrounding obstacles. Similar collision-avoidance techniques have been proposed 
%in~\cite{liu2017planning, schoels2020nmpc} to compute convex safety polytopes for drones using obstacle inflation and point mass representation. 
in~\cite{liu2017planning, schoels2020ciao} to compute convex safety envelopes for drones and linear systems, respectively.
These previously presented techniques can robustly be used for vehicles whose body can be reasonably described by a single bounding circle, i.e., not for elongated vehicles such as ships that are considered in this work. 
Compared to using individual bounding regions for each obstacle, as was done in~\cite{bergman2019bimproved, bergman2020optimization}, the number of obstacle avoidance constraints using the proposed safety-envelope approach scales significantly better in the prediction horizon with an increased number of obstacles. 
Additionally, as the ship's body can be accurately approximated as a convex polytope~\cite{martinsen2019autonomous}, a sufficient condition for collision avoidance at each time instance is to enforce that each vertex of the ship's body is within the computed convex safety envelope. %their corresponding spatial convex polytope.

%In particular, efficient handling of large maps of complex environments with a geometry that is able to capture the one of, e.g., harbors is developed for the second opti. 

The remainder of the paper is organized as follows. In Section~\ref{sec:probform}, the dynamic ship model and the motion planning problem are presented. 
The lattice-based motion planner is presented in Section~\ref{sec:motionplanning}. 
In Section~\ref{sec:control}, the collision-avoidance algorithm and the optimization-based receding-horizon improvement step are presented. 
Simulation results from autonomous docking in a model of the Cape Town harbor are presented in Section~\ref{sec:Res} and the paper is concluded in Section~\ref{sec:conclusions} with a summary of the contributions and a discussion of future research directions.

\section{Ship modeling and problem formulation}\label{sec:probform}
In this section, the dynamic ship model used for motion planning is briefly described and the general motion planning problem is posed as an optimal control problem (OCP). The notation is adopted from~\cite{fossen2011handbook} and for more details regarding ship modeling, the reader is referred to, e.g.,~\cite{fossen2011handbook,perez2006ship}.

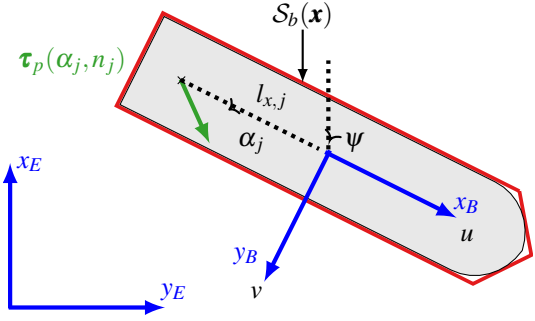
\begin{figure}[!t]
	\centering
	\setlength\figureheight{0.3\textwidth}
	\setlength\figurewidth{0.4\textwidth}
	% This file was created by matlab2tikz.
%
%The latest updates can be retrieved from
%  http://www.mathworks.com/matlabcentral/fileexchange/22022-matlab2tikz-matlab2tikz
%where you can also make suggestions and rate matlab2tikz.
%
\definecolor{mycolor1}{rgb}{0.89020,0.10196,0.10980}%
\definecolor{mycolor2}{rgb}{0.2000,0.6275,0.1725}%

\begin{tikzpicture}

\begin{axis}[%
width=\figurewidth,
height=\figureheight,
at={(0\figurewidth,0\figureheight)},
scale only axis,
xmin=-63.2602939004499,
xmax=39.8511971894683,
ymin=-12.8581762917281,
ymax=12.8786287775784,
axis background/.style={fill=white},
axis x line*=bottom,
axis y line*=left,
axis line style={draw=none},
tick style={draw=none},
ytick=\empty,
xtick=\empty,
axis equal,
]
\addplot [color=mycolor1, line width=1.4pt, forget plot]
table[row sep=crcr]{%
	37.2193889287494	-7.31755117800078\\
	39.5511971894683	-19.7755985947342\\
	28.1856742996503	-25.3849804361991\\
	-41.2602939004499	9.33800366385103\\
	-32.2265792713508	27.4054329220493\\
	37.21938892874940	-7.31755117800078\\
};

\addplot[area legend, draw=black, fill=white!91!black, forget plot]
table[row sep=crcr] {%
	x	y\\
	32.5236461759999	-5.52869679600094\\
	34.1261166902488	-6.52560441982046\\
	35.5233590554552	-7.80529729666178\\
	36.6729187693914	-9.32889260601816\\
	37.5398670418544	-11.050096675391\\
	38.0978620882828	-12.9166115901039\\
	38.32994951155	-14.8717242389754\\
	38.2290774527414	-16.8560295135041\\
	37.7983108583204	-18.8092353020304\\
	37.7263747566166	-18.9531075054381\\
	36.4222700824481	-20.4696442540907\\
	34.8953490981103	-21.7409250658547\\
	33.1920065250528	-22.7283227165639\\
	31.3639975654254	-23.401835628249\\
	29.4668653464494	-23.7409994519022\\
	27.5582532706025	-23.7355088663671\\
	25.6961535500056	-23.3855307003068\\
	23.9371451424007	-22.7016988631993\\
	23.9371451424007	-22.7016988631993\\
	-40.5894735072	9.56161046160101\\
	-32.0029724736008	26.7346125287994\\
	32.5236461759999	-5.52869679600095\\
	32.5236461759999	-5.52869679600094\\
}--cycle;
\addplot [color=black, draw=none, mark=x, mark options={solid, black}, forget plot]
table[row sep=crcr]{%
	-28.6216701119973	14.3108350559987\\
};

\node[] (source) at (axis cs:-29.6,16.3){};
\node (destination) at (axis cs:-22.3,0.15){};
\draw[mycolor2,->, ultra thick, >=latex](source)--(destination);
\node[] (source) at (axis cs:-30.6,15.3){};
\node (destination) at (axis cs:-0,0){};
\draw[dotted, ultra thick, >=latex](source)--(destination);
 \draw[  thick, ->] (axis cs:-19,9) arc [radius=12,start angle=-26.56,end angle=-60];

\node[] (source) at (axis cs:-64.1,-30){};
\node (destination) at (axis cs:-30,-30){};
\draw[blue,->, ultra thick, >=latex](source)--(destination);
\node[] (source) at (axis cs:-62,-32){};
\node (destination) at (axis cs:-62,0){};
\draw[blue,->, ultra thick, >=latex](source)--(destination);
\node[blue] at (axis cs: -30,-27) { $y_E$};
\node[blue] at (axis cs: -58,-2) { $x_E$};
\node[black] at (axis cs: -15,2) { $\alpha_j$};
\node[black] at (axis cs: -11,11) { $l_{x,j}$};
\node[mycolor2] at (axis cs: -50,18) { $\bm{\tau}_p(\alpha_j, n_{j})$};
\node[] (source) at (axis cs:0,-1){};
\node (destination) at (axis cs:0,20){};
\draw[dotted, ultra thick, >=latex](source)--(destination);
\node[] (source2) at (axis cs:-26.83/11,13.42/11){};
\node (destination2) at (axis cs:26.83,-13.42){};
\draw[blue, ->, ultra thick, >=latex](source2)--(destination2);
 \draw[  thick, ->] (axis cs:0,3) arc [radius=3,start angle=90,end angle=-26.56];
\node[] (source) at (axis cs:13.42/11,26.83/11){};
\node (destination) at (axis cs:-13.42,-26.83){};
\draw[blue,->, ultra thick, >=latex](source)--(destination);
\node[] at (axis cs: 5,2.2) { $\psi$};
\node[blue] at (axis cs: 27,-10) { $x_B$};
\node[blue] at (axis cs: -16,-20) { $y_B$};
\node[] at (axis cs: 27,-16) { $u$};
\node[] at (axis cs: -14,-27) { $v$};

\node[] (source) at (axis cs:-5,26){};
\node (destination) at (axis cs:-5,12){};
\draw[->,thick, >=latex](source)--(destination);
\node[] at (axis cs: -5,27) { $\sship(\xbm)$};
\end{axis}
\end{tikzpicture}%
 \vspace*{-1.2em}
	\caption{ \label{fig:ship} \small Definition of Earth-fixed ($x_E,y_E$) and body-fixed ($x_B, y_B$) coordinate systems, body-fixed velocities ($u,v$) and the convex polytope $\sship(\xbm)$ representing the bounding region of the ship's body at $\xbm$. Furthermore, the parameters for the $j^{\mathrm{th}}$ thruster are illustrated. \vspace*{-0.8em}   }
\end{figure}

The ship is assumed to move on the ocean surface and it is thus sufficient to consider the horizontal 3 degrees of freedom motion. The motion of the ship is described using two coordinate systems: an Earth-fixed (inertial) system and a body-fixed system that is located at the center of mass of the ship. The Earth-fixed generalized position is given by $\bm{\eta} = \begin{bmatrix}
x  &y  &\psi\end{bmatrix}^\intercal \in  \mathbf{SE}(2)$~\cite{lavalle2006planning} and the body-fixed generalized velocity vector is represented by $\bm{\nu} = \begin{bmatrix}u  &v  &r\end{bmatrix}^\intercal \in \mathbf{R}^3$, see Fig.~\ref{fig:ship}. The generalized velocity and position is related through the kinematics
\begin{equation}\label{eq:kinematics}
\dot{\bm{\eta}} = \bm{R}(\psi)\bm{\nu}, 
\end{equation}
where
\begin{equation}
\bm{R}(\psi) = \begin{bmatrix}
\cos \psi & -\sin \psi & 0 \\ \sin \psi & \cos \psi & 0 \\ 0 & 0 & 1 \end{bmatrix}
\end{equation}
is the rotation matrix. The kinetics, i.e., the motion induced by forces acting on the ship, is derived using rigid-body mechanics and theory of hydrodynamics, see \cite{fossen2011handbook} for details. 
A model of the kinetics is given by
\begin{equation} \label{eq:dynamics}
\bm{M}\dot{\bm{\nu}} + \bm{C}(\bm{\nu}) + \bm{D}(\bm{\nu}) = \bm{\tau},
\end{equation}
where $\bm{M}$ is the total inertia matrix, $\bm{C}(\bm{\nu})$ corresponds to Coriolis and centripetal forces, $\bm{D}(\bm{\nu})$ describes the damping and $\bm{\tau}$ describes the forces acting on the ship. Simplified models for the thrusters are used where it is assumed that forces induced by the propeller and rudder can be separated. The total forces can be written as $\taubm = \sum_{j=1}^{N_t} \taubm_j$, where $N_t$ is the number of thrusters. The generalized forces of the $j^\text{th}$ thruster is assumed to be
\begin{equation} \label{eq:thrust_alloc}
\bm{\tau}_j =  \underbrace{\begin{bmatrix}1&0& l_{y,j}\\0&1&l_{x,j}\end{bmatrix}^{\sf T}}_{\coloneqq \bm{T}} \left(\theta_{p,j}\bm{\tau}_p({\alpha}_j,{n}_j) + \theta_{r,j}\taubm_r(\alpha_j, \nubm)\right),
\end{equation}
where $l_{x,j}$ and $l_{y,j}$ describe the position of the thruster in the body-fixed coordinate system, $\theta_{p,j}$ is the propeller gain, $\theta_{r,j}$ is the rudder gain, $\alpha_j$ is the thruster angle and $n_j$ is the propeller velocity. The propeller force is assumed to be quadratic with respect to the propeller velocity:
\begin{equation}
\bm{\tau}_{p}(\alpha_j,n_j) = 
\begin{bmatrix}
\cos \alpha_j \\
\sin \alpha_j
\end{bmatrix}|n_j|n_j, \quad j = 1, \ldots, N_t.
\end{equation}
The rudder force is assumed to be proportional to the square of the velocity $V_j$ of the incoming water flow which is dependent on the thruster position and the ship velocity. The rudder force at the $j^{\mathrm{th}}$ thruster is (see, e.g.,~\cite{du2018motion} for details):
\begin{equation}\label{eq:rudderforce}
\bm{\tau}_{r}(\alpha_j,\bm{\nu}) = 
\begin{bmatrix}
-\sin \alpha_j \\
\cos \alpha_j
\end{bmatrix}F_N(V_j,\alpha_j),  \quad j = 1, \ldots, N_t.
\end{equation}
This model is consistent with \cite{Lewandowski} where the lift force is proportional to $V_j^2\alpha_j$ and the drag force is proportional to $V_j^2\alpha_j^2$ for small angles of attack.
To account for dynamics in the thrusters, the derivatives of the thruster angles and velocities are considered as control inputs to the system. A compact notation of the state vector is thus defined as $\bm{x} = \begin{bmatrix}\bm{\eta}^\intercal \; \; \bm{\nu}^\intercal \; \; \alphabm^\intercal \; \; \nbm^\intercal \end{bmatrix}^\intercal$
and the control-input vector is \mbox{$\bm{u} = \begin{bmatrix}\dot{\alphabm}^\intercal &\dot{\nbm}^\intercal \end{bmatrix}^\intercal$}, where \mbox{$\alphabm = \begin{bmatrix}
	\alpha_1  &\ldots  &\alpha_{N_t}\end{bmatrix}^\intercal$} and \mbox{$\nbm = \begin{bmatrix}
	n_1  &\ldots  &n_{N_t}\end{bmatrix}^\intercal$}. 
Then by combining \eqref{eq:kinematics}-\eqref{eq:rudderforce}, the ship dynamics can be written as:

\begin{align}
\dot{\bm{x}}(t) = \bm{f}_\text{dyn}({\bm{x}}(t),{\bm{u}}(t)).
\end{align}
The feasible set for the state space is given by 
\begin{align}
\mathcalOld{X} = \{ \xbm \; | \; ||\begin{bmatrix}
u &v\end{bmatrix}||_2  \leq v_{\mathrm{max}}, \; |n_j| \leq \bar{n}, \quad j = 1, \ldots, N_t \},
\end{align} where $v_{\mathrm{max}}$ is the maximum allowed magnitude of the ship's velocity and $\bar{n}$ is the maximum speed for the thrusters. The control inputs are restricted to
\begin{align}
\mathcalOld{U} = \{\ubm \; | \; |\dot{\alpha}_j| \leq \bar{\dot{\alpha}}_j, \quad |\dot{n}_j| \leq \bar{\dot{n}}_j, \quad j = 1 ,\ldots, N_t \}.
\end{align}  
The ship is assumed to operate in an environment with only static obstacles $\mathcalOld X_{\mathrm{obst}}$. The free-space where the vehicle is not in collision with any obstacle is defined as $\mathcalOld X_{\mathrm{free}} = \mathcalOld X \setminus \mathcalOld X_{\mathrm{obst}}$.
Since the free-space $\mathcalOld X_{\mathrm{free}}$ is defined as the complement set of $\mathcalOld X_{\mathrm{obst}}$, it is in general non-convex. 

The motion planning problem is defined as follows: compute a feasible and collision-free state and control input trajectory $(\bm x(t),\bm u(t))$, $t\in[0,T_f]$ that moves the ship from its initial state ${\bm{x}}(0)=\bm{x}_\text{init}\in \mathcalOld X_{\mathrm{free}}$ to a desired goal state ${\bm{x}}(T_f)=\bm{x}_\text{term}\in \mathcalOld X_{\mathrm{free}}$, while minimizing an objective functional $J$. In this work, the cost function used to define the objective functional is separated into two parts, where $\ell(\xbm,\ubm)$ is the first part that is position invariant, which can, e.g., be related to time and energy consumption. The second part is denoted $\ell_{\mathrm{env}}$ and depends on the generalized position of the ship. It can, e.g., represent costs related to safety. In this work, the cost $\ell_{\mathrm{env}}$ is related to the distance to obstacles:
\begin{equation} \label{eq:distcost}
\begin{aligned}
\ell_{\mathrm{env}}(d) &= \begin{cases}
k_d(d_{\mathrm{safe}}-d)^2, \quad 0 < d \leq d_{\mathrm{safe}} \\
0, \quad\quad\quad\quad\quad\hspace{27pt} d>d_{\mathrm{safe}} 
\end{cases},
\end{aligned}
\end{equation} 
where $d$ is the distance from the bounding region of the ship's body to the closest obstacle in $\mathcalOld{X}_{\mathrm{obst}}$, $d_{\mathrm{safe}}$ is a (soft) desired safety distance and $k_d > 0$ is a weighting parameter. The distance $d$ is formally defined as
\mbox{$d = \mathrm{distance}(\sship(\xbm), \mathcalOld{X}_{\mathrm{obst}})$}, where the bounding region of the ship's body $\sship(\xbm)$ is represented as a convex polytope and depends of the state $\xbm$, see Fig.~\ref{fig:ship}. 
Now, it is possible to pose the motion planning problem as the following continuous-time OCP:
\begin{alignat}{2} 	\label{eq:ocp}
&\minimizer_{\ubm(\cdot), T_f} 
&&J = \int_{0}^{T_f} \left[ \ell(\xbm(t), \ubm(t)  ) + \ell_{\mathrm{env}}(d(t)) \right] \mathrm{d}t  \nonumber \\
&\hspace{0.55em}  \subjecttoo \hspace{1em} 
&& \xbm(0) = \xbm_{\mathrm{init}}, \quad \xbm(T_f)= \xbm_{\mathrm{term}},\nonumber \\ 
& \text{} && \dot{\xbm}(t) = \bm{f}_{\mathrm{dyn}}(\xbm(t), \ubm(t)), \\
& \text{} && d(t) = \mathrm{distance}(\sship(\xbm(t)), \mathcalOld{X}_{\mathrm{obst}}), \nonumber   \\
& \text{} && \xbm(t) \in \mathcalOld{X}_{\mathrm{free}}, \; \ubm(t) \in \mathcalOld{U}.\nonumber 
\end{alignat} 
The upcoming sections describe how the proposed framework for solving \eqref{eq:ocp} is implemented by using a combination of a lattice-based motion planner and an optimization-based improvement step.

\section{Lattice-based motion planner} \label{sec:motionplanning}
A lattice-based motion planner computes a suboptimal solution to the continuous OCP in~\eqref{eq:ocp} by transforming it as a graph-search problem. This transformation is done by restricting the controls to a finite subset of the available actions. In this work, the so-called state-lattice methodology is used~\cite{pivtoraiko2009differentially}, where the set of available actions is represented using a motion-primitive set $\mathcalOld{P}$. The construction of this set is performed offline and can be divided into two main steps: state-space discretization and motion-primitive computation. 

\subsection{State-space discretization}
Before the motion primitive set can be computed, a discrete representation of the ship's state space $\mathcalOld{X}_d$ needs to be selected. Since the complexity of the online graph-search problem grows exponentially with the dimension of the search space~\cite{lavalle2006planning}, it is important to carefully select the discretization $\mathcalOld{X}_d$. In this work, the position $(x,y)$ of the ship is discretized onto a uniform grid with resolution $r_p$. The orientation $\psi \in \Psi$, $|\Psi| = 16$, is irregularly discretized as proposed in~\cite{pivtoraiko2009differentially} to enable the computation of short sway and surge motion primitives. The velocity $\nubm_d \in \mathbf{V}$ is also discretized such that short sway and surge maneuvers can be performed, with rotational velocity $r_d = 0$ for all $\nubm_d$. One such discretization, that will be used in Section~\ref{sec:Res}, is to include surge maneuvers with \mbox{$u_d \in \{0, \pm \frac{v_{\mathrm{max}}}{2}, v_{\mathrm{max}}\}$} and sway maneuvers with $v_d \in \{0, \pm \frac{v_{\mathrm{max}}}{2}\}$.

To be able to compute short sway or surge motion primitives for the discretized values of $\nubm_d \in \mathbf{V}$, the corresponding discretized values $\alphabm_d $ and $\nbm_d $ are selected such that the ship obtains a configuration where $\dot{\nubm}_d = 0$. From the dynamics of the ship given in~\eqref{eq:dynamics}, this configuration can be obtained by selecting $\taubm _d $ such that:
\begin{equation} \label{eq:tau_d}
\bm{M}\dot{\nubm}_d = 0 \iff \taubm _d = - \bm{C}(\nubm _d) - \bm{D}(\nubm _d).
\end{equation}  
Finally, to obtain the thruster angles and velocities $\alphabm_d$ and $\nbm_d$ that generates the desired $\taubm _d$ from \eqref{eq:tau_d}, the following feasibility problem obtained from \eqref{eq:thrust_alloc} is solved:
\begin{alignat}{2} 	\label{eq:ulup}
&\find && \hspace*{0.5em} \alphabm_d, \nbm_d  \nonumber \\
&\hspace*{0.2em}\mathrm{s.t.} 
&& \hspace*{0.5em} \taubm_d = \sum_{i=j}^{M}\bm{T}\left(\theta_{p,j}\bm{\tau}_p({\alpha}_{j,d},{n}_{j,d}) + \theta_{r,j}\taubm_r(\alpha_{j,d}, \nubm_d)\right),  
\end{alignat} 
which gives the complete state-space discretization $\mathcalOld{X}_d$. However, it should be stressed that on the trajectory between two discretized states, the ship can take any feasible state in $\mathcalOld X$.
\subsection{Motion primitive computation}
\begin{figure}[t]
	\centering
	\setlength\figureheight{0.2625\textwidth}
	\setlength\figurewidth{0.35\textwidth}
	\input{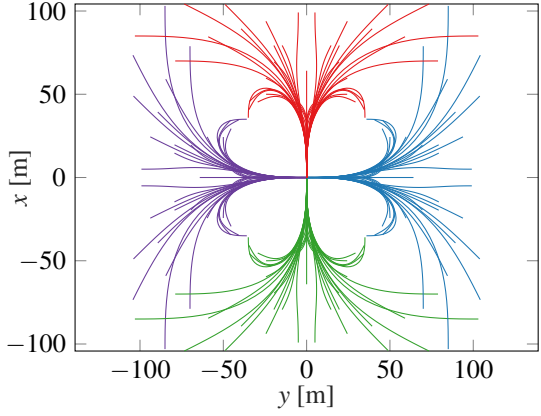}\vspace*{-0.5em}
	\caption{ \label{fig:primitives} \small Illustration of the paths for motion primitives from four different initial headings, $\psi_0 = \frac{k\pi}{2}$,  \mbox{$k= 0,1,2,3$} and \mbox{$\nubm _0 = [1.5433 \mathrm{m/s}, 0 ,0]^\intercal$} to neighboring states in $\mathcalOld{X}_d$. \vspace*{-0.9em}  }
\end{figure}

A motion primitive $\bm{m}$ is defined as a feasible state and control-input trajectory that moves the ship from an initial state $\xbm_0 \in \mathcalOld{X}_d$ to a terminal state $\xbm_f \in \mathcalOld{X}_d$ in time $T$, i.e., 
\begin{equation} \label{eq:motionPrimitive}
\bm{m} =  \big( \xbm(t), \ubm(t)\big)  \in  \mathcalOld{X} \times \mathcalOld{U}, \;   t \in [ 0, T].	\end{equation}
To compute the motion primitives $\bm{m} \in \mathcalOld{P}$, the optimization-based framework in~\cite{bergman2019improved} is used. Since the position-dependent cost ($\ell_{\mathrm{env}}$ in \eqref{eq:ocp}) is unknown offline, the motion primitives are computed by solving OCPs where only the integral of the position invariant cost $\ell_m(\bm{m}) = \int_{0}^{T} \ell(\xbm, \ubm) \mathrm{d}t$ is minimized. 
In this work, the first type of maneuvers are short surge and sway trajectories where the velocity is kept constant. 
The second typ of maneuvers are trajectories to neighboring headings and velocities in $\mathcalOld{X}_d$. 
From each initial state, trajectories to the closest neighbors in $\mathbf{V}$ combined with heading changes to the eight closest neighbors in $\Psi$ are computed. 
To reduce the branching factor of the online graph-search problem, heading changes are limited to the closest adjacent headings from initial states with $u_d = v_{\mathrm{max}}$. 
Moreover, since the ship is orientation invariant, rotational symmetries~\cite{pivtoraiko2005efficient} are exploited to reduce the number of OCPs needed to be solved when $\mathcalOld{P}$ is computed (illustrated in Fig.~\ref{fig:primitives}). 
Together, these maneuvers define the motion primitive set $\mathcalOld{P}$. Since the ship is position invariant, each $\bm m\in\mathcalOld{P}$ can be translated and reused from each discrete position in $\mathcalOld{X}_d$. 
%The resulting maneuvers define the motion primitive set $\mathcalOld{P}$ which can be used to traverse the graph and find a combination of motion that brings the system to all states in $\mathcalOld{X}_d$. 

\begin{figure}[t!]
	\centering
	\setlength\figureheight{0.2625\textwidth}
	\setlength\figurewidth{0.35\textwidth}
	% This file was created by matlab2tikz.
%
%The latest updates can be retrieved from
%  http://www.mathworks.com/matlabcentral/fileexchange/22022-matlab2tikz-matlab2tikz
%where you can also make suggestions and rate matlab2tikz.
%
\definecolor{mycolor1}{rgb}{0.65098,0.80784,0.89020}%
%\definecolor{mycolor2}{rgb}{0.20000,0.62745,0.17255}%
\definecolor{mycolor2}{RGB}{227,26,28}
\begin{tikzpicture}

\begin{axis}[%
width=\figurewidth,
height=\figureheight,
at={(0\figurewidth,0\figureheight)},
scale only axis,
xmin=-20,
xmax=140,
ymin=-50,
ymax=50,
axis background/.style={fill=white},
xlabel style={font=\color{white!15!black}},
xlabel={$y$ [m]},
ylabel style={font=\color{white!15!black}},
ylabel={$x$ [m]},
%axis x line*=bottom,
%axis y line*=left,
ylabel style={yshift=-0.5cm}, %shifting the y line text
xlabel style={yshift=0.2cm}, %shifting the y line text
xmajorgrids,
ymajorgrids,
axis equal
]

\addplot[area legend, draw=mycolor1, fill=mycolor1, forget plot]
table[row sep=crcr] {%
x	y\\
-1.0	-42\\
-2	-42\\
-3	-42\\
-4	-42\\
-5	-42\\
-6	-42\\
-7	-42\\
-8	-42\\
-9	-42\\
-10	-42\\
-11	-41\\
-11	-40\\
-11	-39\\
-11	-38\\
-11	-37\\
-11	-36\\
-11	-35\\
-11	-34\\
-11	-33\\
-11	-32\\
-11	-31\\
-11	-30\\
-11	-29\\
-11	-28\\
-11	-27\\
-11	-26\\
-11	-25\\
-11	-24\\
-11	-23\\
-11	-22\\
-11	-21\\
-11	-20\\
-11	-19\\
-11	-18\\
-11	-17\\
-11	-16\\
-11	-15\\
-11	-14\\
-11	-13\\
-11	-12\\
-11	-11\\
-11	-10\\
-11	-9\\
-11	-8\\
-11	-7\\
-11	-6\\
-11	-5\\
-11	-4\\
-11	-3\\
-11	-2\\
-11	-1\\
-11	0\\
-11	1\\
-11	2\\
-11	3\\
-11	4\\
-11	5\\
-11	6\\
-11	7\\
-11	8\\
-11	9\\
-11	10\\
-11	11\\
-11	12\\
-11	13\\
-11	14\\
-11	15\\
-11	16\\
-11	17\\
-11	18\\
-11	19\\
-11	20\\
-11	21\\
-11	22\\
-11	23\\
-11	24\\
-11	25\\
-11	26\\
-11	27\\
-11	28\\
-11	29\\
-11	30\\
-11	31\\
-11	32\\
-11	33\\
-11	34\\
-11	35\\
-11	36\\
-10	37\\
-9	38\\
-8	39\\
-7	39\\
-6	40\\
-5	41\\
-4	42\\
-3	42\\
-2	43\\
-1	44\\
0	45\\
1	45\\
2	45\\
3	45\\
4	45\\
5	45\\
6	44\\
7	44\\
8	44\\
9	45\\
10	44\\
11	44\\
12	44\\
13	44\\
14	44\\
15	44\\
16	44\\
17	44\\
18	44\\
19	44\\
20	44\\
21	44\\
22	44\\
23	44\\
24	44\\
25	44\\
26	43\\
27	43\\
28	43\\
29	43\\
30	42\\
31	43\\
32	42\\
33	42\\
34	42\\
35	42\\
36	41\\
37	42\\
38	41\\
39	41\\
40	40\\
41	41\\
42	40\\
43	40\\
44	39\\
45	39\\
46	39\\
47	39\\
48	38\\
49	38\\
50	37\\
51	37\\
52	36\\
53	36\\
54	36\\
55	36\\
56	35\\
57	35\\
58	34\\
59	34\\
60	33\\
61	33\\
62	32\\
63	31\\
64	31\\
65	30\\
66	30\\
67	29\\
68	29\\
69	28\\
70	28\\
71	27\\
72	26\\
73	26\\
74	25\\
75	24\\
76	24\\
77	23\\
78	23\\
79	23\\
80	22\\
81	22\\
82	21\\
83	21\\
84	21\\
85	20\\
86	20\\
87	20\\
88	19\\
89	19\\
90	19\\
91	19\\
92	18\\
93	18\\
94	18\\
95	18\\
96	17\\
97	17\\
98	17\\
99	17\\
100	17\\
101	17\\
102	16\\
103	16\\
104	16\\
105	16\\
106	16\\
107	16\\
108	16\\
109	16\\
110	16\\
111	16\\
112	16\\
113	16\\
114	16\\
115	16\\
116	16\\
117	16\\
118	16\\
119	16\\
120	16\\
121	16\\
122	16\\
123	16\\
124	16\\
125	16\\
126	16\\
127	15\\
128	14\\
129	13\\
129	12\\
130	11\\
131	10\\
132	9\\
132	8\\
133	7\\
134	6\\
135	5\\
134	4\\
133	3\\
132	2\\
132	1\\
131	0\\
130	-1\\
129	-2\\
129	-3\\
128	-4\\
127	-5\\
126	-6\\
125	-6\\
124	-6\\
123	-6\\
122	-6\\
121	-6\\
120	-6\\
119	-6\\
118	-6\\
117	-6\\
116	-6\\
115	-6\\
114	-6\\
113	-6\\
112	-6\\
111	-6\\
110	-6\\
109	-6\\
108	-6\\
107	-6\\
106	-6\\
105	-6\\
104	-6\\
103	-6\\
102	-6\\
101	-6\\
100	-6\\
99	-6\\
98	-6\\
97	-6\\
96	-6\\
95	-6\\
94	-6\\
93	-6\\
92	-6\\
91	-6\\
90	-6\\
89	-6\\
88	-6\\
87	-6\\
86	-6\\
85	-6\\
84	-6\\
83	-6\\
82	-6\\
81	-6\\
80	-6\\
79	-6\\
78	-6\\
77	-6\\
76	-6\\
75	-6\\
74	-6\\
73	-6\\
72	-6\\
71	-6\\
70	-6\\
69	-6\\
68	-6\\
67	-6\\
66	-6\\
65	-6\\
64	-6\\
63	-6\\
62	-6\\
61	-7\\
60	-7\\
59	-7\\
58	-7\\
57	-7\\
56	-7\\
55	-7\\
54	-7\\
53	-7\\
52	-8\\
51	-8\\
50	-8\\
49	-8\\
48	-8\\
47	-9\\
46	-9\\
45	-9\\
44	-9\\
43	-10\\
42	-10\\
41	-10\\
40	-10\\
39	-11\\
38	-11\\
37	-12\\
36	-12\\
35	-13\\
34	-13\\
33	-14\\
32	-14\\
31	-15\\
30	-16\\
29	-16\\
28	-17\\
27	-18\\
27	-19\\
26	-20\\
25	-21\\
25	-22\\
24	-23\\
24	-24\\
23	-25\\
23	-26\\
22	-27\\
22	-28\\
22	-29\\
21	-30\\
21	-31\\
21	-32\\
21	-33\\
21	-34\\
20	-35\\
20	-36\\
20	-37\\
20	-38\\
20	-39\\
19	-40\\
19	-41\\
19	-42\\
18	-43\\
17	-43\\
16	-42\\
15	-42\\
14	-42\\
13	-42\\
12	-42\\
11	-42\\
10	-42\\
9	-42\\
8	-42\\
7	-42\\
6	-42\\
5	-42\\
4	-42\\
3	-42\\
2	-42\\
1	-42\\
0	-42\\
}--cycle;
\addplot [color=mycolor2, line width=1.8pt, forget plot]
  table[row sep=crcr]{%
-10.1	36.5625418628434\\
0	44.2195827535749\\
10.1	36.5625418628434\\
10.1	-41.080410966513\\
-10.1	-41.080410966513\\
-10.1	36.5625418628434\\
};

\addplot[area legend, draw=black, fill=white!91!black, forget plot]
table[row sep=crcr] {%
x	y\\
-9.6	31.5625418628434\\
-9.42498531470064	33.4416657063929\\
-8.90525899145366	35.2636933228384\\
-8.05661265159672	36.9732633250452\\
-6.90483198473064	38.5184312935799\\
-5.48491326358891	39.8522480816118\\
-3.84	40.9341863409803\\
-2.0200720511117	41.7313719251613\\
-0.0804270067230374	42.2195827535749\\
0.0804270067230374	42.2195827535749\\
2.0200720511117	41.7313719251613\\
3.84	40.9341863409803\\
5.48491326358891	39.8522480816118\\
6.90483198473064	38.5184312935799\\
8.05661265159672	36.9732633250452\\
8.90525899145366	35.2636933228384\\
9.42498531470064	33.4416657063929\\
9.6	31.5625418628434\\
9.6	31.5625418628434\\
9.6	-40.580410966513\\
-9.6	-40.580410966513\\
-9.6	31.5625418628434\\
-9.6	31.5625418628434\\
}--cycle;
\addplot [color=mycolor2, line width=1.8pt, forget plot]
  table[row sep=crcr]{%
42.7533838575576	35.3182276092109\\
55.4267241791493	35.1547641715505\\
58.6916760163262	22.9081170915025\\
10.9908034651612	-38.3540644067039\\
-4.94748869360734	-25.9439538889955\\
42.7533838575576	35.3182276092109\\
};

\addplot[area legend, draw=black, fill=white!91!black, forget plot]
table[row sep=crcr] {%
x	y\\
40.0760865056191	31.0659248343249\\
41.3686394640986	32.4410769746324\\
42.8981005566193	33.559401101454\\
44.61799788183	34.3869175020959\\
46.4760732298137	34.8984825094676\\
48.4158699213539	35.078552479685\\
50.3784482155625	34.9216560778227\\
52.3041761640823	34.4325605216576\\
54.134541497663	33.6261267321696\\
54.2614592333667	33.5273041537496\\
55.4919482960013	31.9504498871992\\
56.4381553333518	30.2033566334662\\
57.071330323392	28.3391089369936\\
57.372234564452	26.4143510206839\\
57.3317252334319	24.48756567985\\
56.9510331858559	22.6172973120722\\
56.2417255569503	20.8603730753044\\
55.2253543000922	19.2701762234337\\
55.2253543000922	19.2701762234337\\
10.9034722364221	-37.6523712711476\\
-4.24579555805099	-25.8566226602564\\
40.0760865056191	31.0659248343249\\
40.0760865056191	31.0659248343249\\
}--cycle;
\addplot [color=mycolor2, line width=1.8pt, forget plot]
  table[row sep=crcr]{%
126.562538562213	15.1000119484515\\
134.219582753573	5.00001445073397\\
126.56254516347	-5.09998805154738\\
48.9195923341177	-5.10001342486811\\
48.9195857328607	15.0999865751308\\
126.562538562213	15.1000119484515\\
};

\addplot[area legend, draw=black, fill=white!91!black, forget plot]
table[row sep=crcr] {%
x	y\\
121.562538725611	14.6000103144771\\
123.441662626354	14.4249962432658\\
125.263690412643	13.9052705154482\\
126.973260692183	13.0566247342701\\
128.518429037114	11.904844572357\\
129.852246289168	10.4849262870999\\
130.934185086086	8.84001337708299\\
131.73137126501	7.02008568871096\\
132.21958272729	5.08044080386721\\
132.219582779856	4.91958679042114\\
131.731372585308	2.97994158648778\\
130.934187595871	1.1600133770834\\
129.852249874051	-0.484900240077329\\
128.518433550042	-1.9048193971035\\
126.973265957903	-3.05660056892253\\
125.26369623303	-3.90524746745817\\
123.441668786428	-4.42497438613445\\
121.562545000073	-4.59998968552188\\
121.562545000073	-4.59998968552188\\
49.4195921707202	-4.60001326147069\\
49.4195858962582	14.5999867385283\\
121.562538725611	14.6000103144771\\
121.562538725611	14.6000103144771\\
}--cycle;
\node[] at (axis cs: 0,0) {$\etabm_0$};
\node[] at (axis cs: 28.2599,0.2644) {};
\node[] at (axis cs: 90,5) {$\etabm_f$};
\node[] (source) at (axis cs:110,40){};
\node (destination) at (axis cs:110,15){};
\draw[->, thick, >=latex](source)--(destination);
\node[] at (axis cs: 110,45) {\small $\sship(\bm{x}_f)$};
\end{axis}
\end{tikzpicture}%
%\vspace*{-0.5em}
	\caption{ \label{fig:swath} \small The swath (blue area) of a motion primitive that moves the ship from $\etabm _0 = [0,0,0]^\intercal$, $\nubm _0 = [0,1.5433 \mathrm{m/s},0]^\intercal$, to  \mbox{$\etabm _f = [5\mathrm{m},90\mathrm{m},\pi/2]^\intercal$}, $\nubm _f = [1.5433\mathrm{m/s},0,0]^\intercal$. \vspace*{-0.9em}  }
\end{figure}
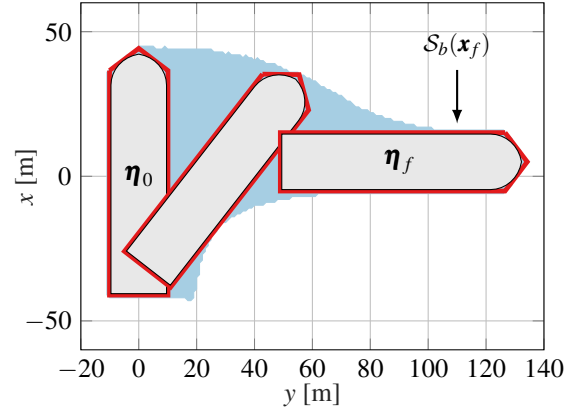
\vspace*{-0.5em}
\subsection{Collision avoidance and environmental cost}\label{sec:swath_comp}
The environment used within the lattice-based motion planner is selected to be represented using a cost map~\cite{ferguson2008efficient}. The cost map is both used to ensure collision avoidance and to obtain an approximation of the position-dependent cost $\ell_{\mathrm{env}}$ in \eqref{eq:ocp}. For efficient online collision checking, we precompute \textit{swaths}~\cite{pivtoraiko2009differentially} for each motion primitive in $\mathcalOld{P}$. A swath is represented as the set of cells that the volume of the ship is occupying during an execution of a motion primitive. In this work, a convex polytope $\sship$ is used to represent the volume of the ship as shown in Fig.~\ref{fig:ship}. One example of a swath for a right-turn maneuver is illustrated in Fig.~\ref{fig:swath}.
Furthermore, it is also possible to approximate the position-dependent cost for $\bm{m}$ applied from $\xbm$ as
\begin{equation}
\ell_{m,\mathrm{env}}(\xbm, \bm{m}, \mathcalOld{X}_{\mathrm{obst}}) =\int_{0}^{T} \ell_{\mathrm{env}}(d(t)) \mathrm{d}t,
\end{equation}
where $d(t) = \mathrm{distance}(\sship(\xbm, \bm{m}(t)), \mathcalOld{X}_{\mathrm{obst}})$. Since only static obstacles are considered in this work, it is possible to precompute an approximation of this function. This approximation is then used as the cost map during planning.  
\subsection{Resulting graph-search problem}
Now, the resulting graph-search problem can be posed as a discrete OCP in the form:
\begin{subequations}
	\label{eq:lattice}
	\begin{alignat}{2}
	&\minimizer_{\{ m_k \}_{k=0}^{M-1},\hspace{0.5ex}M}
	&&J_m = \sum_{k=0}^{M-1} \ell _m(\bm{m}_{k}) + \ell _{m,\mathrm{env}} (\xbm_k, \bm{m}_k, \mathcalOld{X}_\mathrm{obst}) \label{eq:objLattice}  \\
	&\subjecttoo \hspace{1em} 
	&& \xbm_0 = \xbm_{\mathrm{init}}, \quad \xbm_{M} = \xbm_{\mathrm{term}},  \\ \label{eq:primStateTrans}
	& \text{} && \xbm_{k+1} = \bm{f}_{m}(\xbm_k, \bm{m}_{k}), \\ \label{eq:primInSet}
	& \text{} && \bm{m}_{k}  \in \mathcalOld{P}(\xbm_k), \\ \label{eq:trajInFreeSpace}
	& \text{} && c(\xbm_k, \bm{m}_k) \notin \mathcalOld{X}_{\mathrm{obst}},
	\end{alignat} 
\end{subequations}
where the decision variables are the motion-primitive sequence $\{ m_k \}_{k=0}^{M-1}$ and its length $M$. 
The state transition equation in \eqref{eq:primStateTrans} outputs the successor state $\xbm_{k+1}$ when $\bm{m}_k$ is applied from $\xbm _k$, and \eqref{eq:primInSet} restricts the choice of $\bm{m}_k$ to the set of applicable motion primitives $\mathcalOld{P}(\xbm_k)$ from state $\xbm _k$. Finally, the constraint in \eqref{eq:trajInFreeSpace} ensures that the ship does not collide with obstacles when $\bm{m}_k$ is applied from $\xbm_k$. 

The problem in \eqref{eq:lattice} can be solved using graph-search algorithms~\cite{pivtoraiko2009differentially}. For a faster search, a precomputed free-space heuristic look-up table (HLUT)~\cite{knepper2006high} is used as heuristic function. 
It is computed by exploiting rotational symmetries and solving~\eqref{eq:lattice} in free-space to all neighboring states within a predefined distance. Similar to the computation of $\mathcalOld{P}$, rotational symmetries are exploited to reduce the size of the HLUT.  
The output from solving \eqref{eq:lattice} is hereafter referred to as the \emph{nominal trajectory} ($\bar{\xbm}(\cdot), \; \bar{\ubm}(\cdot)$) which is used to warm-start the optimization-based improvement~step.

\section{Optimization-based improvement} \label{sec:control} 
In this section, it is shown how the nominal trajectory computed by the lattice-based planner can be used as warm-start to a receding-horizon optimization-based improvement step. The aim is to remove the suboptimality introduced when solving \eqref{eq:lattice} and to efficiently find a locally optimal solution to the original motion planning problem in~\eqref{eq:ocp}.  

\subsection{Collision avoidance and environmental cost}
One important aspect for reducing the complexity of the optimization-based improvement step is to use an efficient representation of the available free space. One approach is to use individual bounding regions of each obstacle~\cite{lavalle2006planning}, but the complexity using this approach scales poorly when the number of obstacles increases. 
Instead, a more computationally efficient approach is used where so-called safety envelopes~\cite{liu2017planning,schoels2020ciao} are computed, which represent local inner-approximations of the ship's obstacle-avoidance constraints in~\eqref{eq:ocp}. However, previous methods for computing these safety envelopes have only shown robustness for vehicles that can be described by a single bounding circle, i.e., not for ships that is considered in this work.   

By using direct methods for optimal control, the improvement OCP is constructed by discretizing the continuous-time motion planning problem~\eqref{eq:ocp}, e.g., using direct multiple-shooting with $N+1$ discretization points~\cite{diehl2006fast}. At each discrete point $i$, $i = 0, \ldots, N$, a local spatial constraint will be computed for collision avoidance that is represented by a convex polytope $\sspatial^i$ which is defined as
\begin{equation} \label{eq:sspatial}
\sspatial^i = \{ \bm{p} \in \mathbf{R}^2 \; | \; \bm{A}_i \bm{p} \leq \bm{b}_i \},  
\end{equation}
where $\bm{A}_i = [\bm{a}_{i,1} \; \ldots \; \bm{a}_{i,K_i} ]^\intercal \in \mathbf{R}^{K_i \times 2}$, $\bm{a}_{i,j} \in \mathbf{R}^2, \; ||\bm{a}_{i,j}||_2 = 1,$ for $j = 1 ,\ldots ,K_i$, $\bm{b}_i \in \mathbf{R}^{K_i}$ and $K_i$ is the number of half spaces that defines $\sspatial^i$. These polytopes are computed using Algorithm~\ref{alg:ppas}, where each polytope $\sspatial^i$ is temporarily transformed to a directed cyclic graph \mbox{$\gspatial^i = <\vspatial^i, \espatial^i>$}, where the vertices $\bm{p}_k \in \vspatial^i$ define the convex hull of $\sspatial^i$ and each edge $e_k \in \espatial^i$ represents a vector from $\bm{p}_k$ to $\bm{p}_{k+1}$. Each vertex $\bm{p}_k$ is associated with an expansion direction $\bar{\bm{g}}_k$ which is selected such that the area spanned by the graph increases after each expansion. One such example is to select the directions $\bar{\bm{g}}_k$ to be aligned with the vectors pointing from the center of the ship's body to each of the initial vertices.

\begin{figure}
	\centering
	\setlength\figureheight{0.2625\textwidth}
	\setlength\figurewidth{0.35\textwidth}
	\input{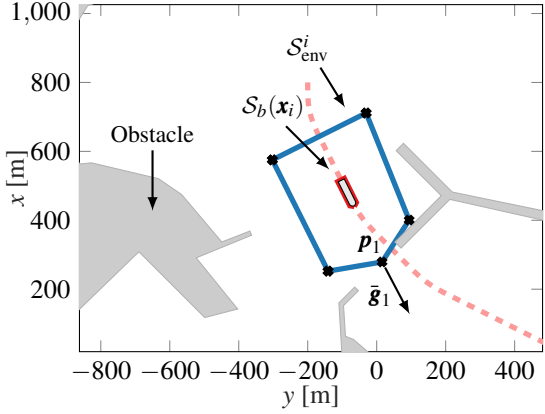} \vspace*{-0.5em}
	\caption{ \label{fig:bb} \small Illustration of one resulting convex bounding box ($\sspatial^i$) computed from a state on the nominal trajectory (dashed). Collision avoidance is ensured by keeping the vertices of the ship $\sship(\xbm_i)$ (red) within the local spatial constraints $\sspatial^i$ (blue).  \vspace*{-0.8em}}
\end{figure}

The inputs to the algorithm are the nominal state from the lattice-based planner at point $i$, i.e., $\etabm_i$ and $\nubm_i$, the polytope $\sship$ representing the ship's body and the edge set $\eworld$ that represents surrounding obstacles. The graph is initialized using $\sship$ at $\etabm_i$ (Line~\ref{line:g_init}). Then, the initial candidate vertex $\bm{p}_c$ is selected as the one whose expansion direction $\bar{\bm{g}}$, see Fig.~\ref{fig:bb}, is closest to the ship's velocity vector $\nubm_i$. 
The expansion of the graph $\gspatial^i$ continues as long as at least one vertex is expandable.
A vertex is no longer considered expandable if its corresponding step length $\Delta_l$ is below a minimum value or if a maximum expansion distance is reached. If a candidate vertex is expandable, it is expanded in its expansion direction $\bar{\bm{g}}$ using step size $\Delta_l$ (Line~\ref{line:expand}). 

After a vertex expansion, it is first checked that the area spanned by the graph is convex. If this is not the case, the algorithm will consume the limiting vertex if the cardinality of $\vspatial^i$ is above a certain threshold (Line~\ref{line:thres}), and the reduced vertex set does not intersect any obstacle (Line~\ref{line:remove_vertex}). Otherwise, the step is rejected, $\Delta_l$ is halved and the expansion continues from its child (Line~\ref{line:continue}). The next step is to calculate the number of intersections $n_{\text{lines}}$ between $\espatial^i$ and $\eworld$ (Line~\ref{line:check_inters}). Depending on the value of $n_{\text{lines}}$, three cases can occur:
\begin{enumerate}
	\item \textit{Exactly one edge in $\eworld$ intersected} (Line~\ref{line:one_intesection}). The current vertex is first projected to feasibility by reducing the step size such that $\bm{p}_c$ lies on the intersected edge. Then, it is duplicated and inserted in the set of vertices. The expansion directions for $\bm{p}_c$ and $\bm{p}_{c,\mathrm{new}}$ are modified such that they grow in each direction parallel to the edge that was intersected.
	\item \textit{More than one edge in $\eworld$ intersected} (Line~\ref{line:more_inters}). The step is rejected and the step size $\Delta_l$ is halved.
	\item Otherwise, the step is accepted and the expansion continues with the next vertex in $\vspatial^i$ (Line~\ref{line:next_vertex}	).
\end{enumerate}
Fig.~\ref{fig:bb} illustrates how Algorithm~\ref{alg:ppas} computes a local spatial polytope $\sspatial^i$ in an area with obstacles using the nominal trajectory computed by the lattice-based planner. 
\begin{algorithm}[t]
	\caption{Computation of local spatial constraints}
	\label{alg:ppas}
	\begin{algorithmic}[1]
		\State \textbf{Input}: $\bar{\etabm}_i$, $\bar{\nubm}_i$, $\sship$, $\eworld$.
		\State $<\vspatial^i, \espatial^i > $ = initialize\_graph($\sship, \bar{\etabm}_i$)  \label{line:g_init}
		\State $\bm{p}_c \leftarrow$  select\_initial($\vspatial^i$, $\bar{\nubm}_i$)
		\While{any($\bm{p}_k$.expandable(), $\forall \bm{p}_k \in \vspatial$) }
		\If{$\bm{p}_c$.expandable()}
		\State $\Delta_l \leftarrow \bm{p}_c$.get\_step\_length()
		\State $\bar{\bm{g}} \leftarrow \bm{p}_c.$get\_direction() 
		\State $\bm{p}_c \leftarrow \bm{p}_c + \Delta_l\bar{\bm{g}}$ \label{line:expand}
		
		\If{span($\vspatial^i$) \textbf{not} convex}
		\State success = \textbf{False}
		\If{$| \vspatial^i | > $ $n_{\mathrm{card}}$}  \label{line:thres}
		\State success $\leftarrow$ $\vspatial^i$.remove\_vertex() \label{line:remove_vertex}
		\EndIf
		\If{\textbf{not} success}
		\State $\bm{p}_c \leftarrow \bm{p}_c - \Delta_l\bar{\bm{g}}$ 
		\State $\bm{p}_c$.set\_step\_length($\Delta_l$/2)
		\State \textbf{continue} \label{line:continue}
		\EndIf
		\EndIf
		\State $n_{\text{lines}} \leftarrow$  check\_intersections($\espatial^i$, $\eworld$) \label{line:check_inters}
		\If{$n_{\text{lines}} = 1$} \label{line:one_intesection}
		\State $\bm{p}_c \leftarrow$ project\_to\_feasibility($\bm{p}_c$)
		\State $\bm{p}_{c,\mathrm{new}} \leftarrow \bm{p}_c $
		\State update\_directions($\bm{p}_{c,\mathrm{new}}$, $\bm{p}_c$) 
		\State $\vspatial$.insert\_vertex($\bm{p}_{c,\mathrm{new}}$)
		\ElsIf{$n_{\text{lines}} > 1$} \label{line:more_inters}
		\State $\bm{p}_c \leftarrow \bm{p}_c - \Delta_l\bar{\bm{g}}$ 
		\State $\bm{p}_c$.set\_step\_length($\Delta_l$/2)
		\EndIf
		\EndIf
		\State $\bm{p}_c$ = $\bm{p}_c$.get\_child() \label{line:next_vertex}
		\EndWhile
	\end{algorithmic}
\end{algorithm}

The output from Algorithm~\ref{alg:ppas} is a convex polytope $\sspatial^i$ that is computed for each discrete point $i$ of the improvement step. Since both $\sship(\xbm_i)$ and $\sspatial^i$ are convex, a sufficient condition to ensure collision avoidance for the ship's body is that all vertices $\vship$ of $\sship(\xbm_i)$ lie within $\sspatial^i$~\cite{martinsen2019autonomous}. This condition can be described by the following constraint:
\begin{equation} \label{eq:ocp_coll}
\bm{A}_i\bm{c}_{\text{rot}}(\xbm_i, \bm{p}_j) \leq \bm{b}_i, \quad \forall \bm{p}_j \in \vship,
\end{equation}
where 
\begin{align}
\bm{c}_{\text{rot}}(\xbm_i, \bm{p}_j) =\begin{bmatrix}
\cos \psi_i & -\sin \psi_i \\ \sin \psi_i & \cos \psi_i
\end{bmatrix}\bm{p}_j + \begin{bmatrix}
x_i \\ y_i
\end{bmatrix},
\end{align} and $\bm{A}_i$ and $\bm{b}_i$ are the half-space representation of $\sspatial^i$ defined in~\eqref{eq:sspatial}. Furthermore, we want  to express the position-dependent cost related to the distance to obstacles also in the improvement step. As a consequence, the following modification of \eqref{eq:ocp_coll} is introduced:
\begin{equation} \label{eq:ocp_coll_epsi}
\begin{aligned}
\bm{A}_i\bm{c}_{\text{rot}}(\xbm_i, \bm{p}_j) &\leq \bm{b}_i - \bm{1}_i(d_{\mathrm{safe}}-\epsilon_{d,i}), \quad \forall \bm{p}_j \in \vship,
\end{aligned}
\end{equation}
where a variable $\epsilon_{d,i}$ which satisfies $0 \leq \epsilon_{d,i} \leq d_{\mathrm{safe}}$ is added, which is used to represent the distance to the boundaries of $\sspatial^i$. This modification makes it possible to compute an approximation of $\ell_{\mathrm{env}}$ in \eqref{eq:ocp} since the distance to obstacles is $d_i \geq d_{\mathrm{safe}}-\epsilon_{d,i}$ when \eqref{eq:ocp_coll_epsi} holds. Note that the feasible sets for~\eqref{eq:ocp_coll} and~\eqref{eq:ocp_coll_epsi} are equal, as~\eqref{eq:ocp_coll} is obtained from~\eqref{eq:ocp_coll_epsi} by selecting $\epsilon_{d,i} = d_{\mathrm{safe}}$. %which potentially could lead to feasibility issues.

\subsection{Improvement optimization problem}
%It is here assumed that a nominal trajectory $\bar{\xbm}(\cdot), \; \bar{\ubm}(\cdot)$ exists, which is initially provided by the lattice-based planner. 
The improvement step can now be posed as the following discrete OCP at the current state of the ship defined by $\bm{x}_{\mathrm{cur}}$:
\begin{alignat}{2}
&\minimizer_{\{ \ubm_i, \epsilon_{d,i} \}_{i=0}^{N-1},\hspace{0.5ex}\Delta_t}
&&J_d = \sum_{i=0}^{N-1} \ell _d(\xbm_{i}, \ubm_i, \Delta_t) + \ell_{d,\mathrm{env}}(d_{\mathrm{safe}} - \epsilon_{d,i}, \Delta_t) \nonumber  \\
&\hspace{0.3em} \subjecttoo 
&& \xbm_0 = \xbm_{\mathrm{cur}}, \quad \xbm_{N} = \bar{\xbm}(\xi), \nonumber  \\ \label{eq:ocp2}
& \text{} && \xbm_{i+1} = \bm{f}_{d}(\xbm_i, \ubm_{i}, \Delta_t), \\
& \text{} && 	\bm{A}_i\bm{c}_{\text{rot}}(\xbm_i, \bm{p}_j)  \leq \bm{b}_i - \bm{1}_i(d_{\mathrm{safe}} - \epsilon_{d,i}), \; \forall \bm{p}_j \in \vship \nonumber\\ \label{eq:hehe}
& \text{} && \xbm_i \in \mathcalOld{X} \quad \ubm _i \in \mathcalOld{U}, \quad 0 \leq \epsilon_{d,i} \leq d_{\mathrm{safe}}. \nonumber
\end{alignat} 
Here, $\Delta_t$ is a decision variable that represents the time between two consecutive discrete points. The objective function $J_d$ is the sum of the numerical integrals of $\ell$ and $\ell_{\mathrm{env}}$ in \eqref{eq:ocp} over the time interval $\Delta_t$. Furthermore, the parameter $\xi$ is used to decide where to connect to the nominal trajectory. Hence, it is possible to solve the problem once  with $\xbm_0 = \xbm_{\mathrm{init}}$ and $\xbm_N = \bar{\xbm}(T_f) = \xbm_{\mathrm{term}}$, or repeatedly using a receding-horizon approach with $\xi$ and $\xbm_{\mathrm{cur}}$ updated at each iteration to maintain a desired planning horizon $T=N\Delta_t$~\cite{bergman2020optimization}. In both cases, the nominal trajectory is used both as warm-start, but also to compute the local convex spatial constraints using Algorithm~\ref{alg:ppas}. Since the initial nominal trajectory is feasible, the improvement step is guaranteed to find a solution which is at least as good as the initial one~\cite{bergman2020optimization}. 
As previously mentioned, the spatial constraints in \eqref{eq:ocp2} are only \emph{local} approximations of the free space computed from the nominal trajectory. To avoid unnecessary computational burden when a receding-horizon approach is used, these spatial constraints are only re-computed if any of the constraints in~\eqref{eq:ocp_coll} are close to become active. Otherwise, the constraints from the previous iteration are reused in the next iteration.

Finally, note that the proposed improvement step can be used to transform the optimal open-loop control law computed by the lattice-based planner to a optimal closed-loop control law if it is repeated at a sufficiently high rate. This will make it possible to suppress various disturbances acting on the ship and changes in the environment during trajectory execution. Another possibility is to let changes in the environment trigger replanning and use a trajectory-tracking controller such as~\cite{barlund2020nonlinear} to stabilize the ship around the more or less static trajectory computed by the improvement step.  

\section{Simulation study} \label{sec:Res}
In this section, the proposed framework is evaluated in a simulation study by solving motion planning problems in a model of the Cape Town harbor. The geometry of the harbor is represented using a polygon that has been calculated using data from Google maps. The cost map used in the lattice-based planner is computed with resolution 1 m. The lattice-based planner is implemented in C++ and the improvement step in Python using CasADi~\cite{andersson2018casadi} with IPOPT~\cite{wachter2006implementation} as NLP solver, and discretized using multiple shooting.
\begin{figure}
	\centering
	\setlength\figureheight{0.2625\textwidth}
	\setlength\figurewidth{0.35\textwidth}
	\input{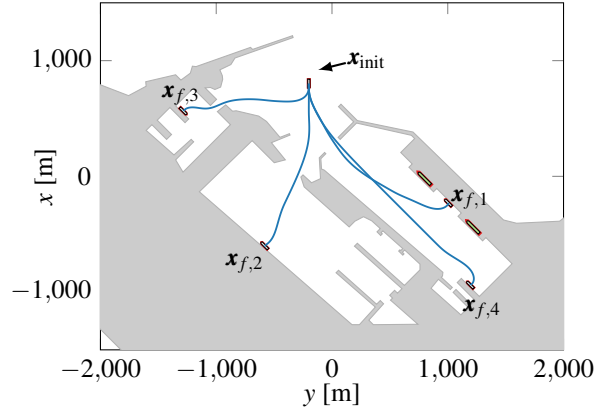} \vspace*{-0.5em}
	\caption{ \label{fig:sc} \small Simulation-study scenario with resulting paths (blue). Planning from an initial state on the fairway to four terminal docking states $\xbm_{f,j}$, $j=1, \ldots, 4$ in the harbor. \vspace*{-0.8em} }
\end{figure}

The simulations are based on a model of a small ship with two thrusters. The model of the ship is based on the \emph{supply vessel} model in the MSS hydro toolbox~\cite{perez2006overview}. The parameters for the two thrusters are:
\begin{equation}
\begin{aligned}
&l_x = \pm 32 \; \mathrm{m}, \quad l_y = 0, \quad \bar{n} = 2 \; \mathrm{RPS}, \\ &\bar{\dot{n}} = 0.08 \; \mathrm{RPS/s}, \quad \bar{\dot{\alpha}} = 7.2 ^\circ\mathrm{/s},	
\end{aligned}
\end{equation}
and the propeller and rudder gains are $\theta_p = 4.2 \cdot 10^{5}$ and $\theta_r = 3.8 \cdot 10^{4}$, respectively. 
Hence, one thruster is mounted in the stern, as shown in Fig.~\ref{fig:ship}, and the other in the bow. The lattice-based planner is discretized with a position resolution of $r_p = 5$ m, and the maximum velocity $v_{\mathrm{max}}$ is set to \mbox{$6$ kn}.
The position-invariant cost function used in both steps of the framework is selected as 
\begin{equation}
\ell(\xbm, \ubm) = 1 + 0.1\nbm^{\sf T}\nbm + 100\dot{\nbm}^{\sf T}\dot{\nbm} + 100\dot{\alphabm}^{\sf T}\dot{\alphabm},
\end{equation}
to trade-off between time, energy and smoothness. 
For the position-dependent cost $\ell_{\mathrm{env}}$, the value of weighting parameter is selected as \mbox{$k_d = 1.5\cdot 10^{-3}$ } in \eqref{eq:distcost} and \eqref{eq:ocp2}, and the safety distance parameter is set to $d_{\mathrm{safe}} = 20 $m which is approximately equal to the width of the ship's body.   

The motion-planning scenario is illustrated in Fig.~\ref{fig:sc}, where four different trajectories are computed from an initial state on the fairway. The terminal states $\xbm_{f,j}$ represent docking states with zero velocity. Since the performance of the proposed motion planner is evaluated in this section, the simulations are performed without any disturbances and the open-loop control law is used for evaluation. 

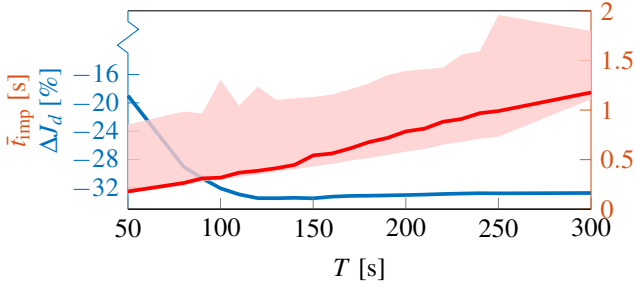
\begin{figure}
	\centering
	\setlength\figureheight{0.15\textwidth}
	\setlength\figurewidth{0.35\textwidth}
	% This file was created by matlab2tikz.
%
%The latest updates can be retrieved from
%  http://www.mathworks.com/matlabcentral/fileexchange/22022-matlab2tikz-matlab2tikz
%where you can also make suggestions and rate matlab2tikz.
%
\definecolor{mycolor1}{rgb}{0.00000,0.44700,0.74100}%
\definecolor{mycolor2}{rgb}{0.85000,0.32500,0.09800}%
\begin{tikzpicture}

\begin{axis}[%
width=\figurewidth,
height=\figureheight,
at={(0\figurewidth,0\figureheight)},
scale only axis,
xmin=50,
xmax=300,
separate axis lines,
every outer y axis line/.append style={mycolor1},
every y tick label/.append style={font=\color{mycolor1}},
every y tick/.append style={mycolor1},
ymin=-35,
ymax=-7,
xmajorticks=false,
ylabel style={yshift=-0.2cm}, %shifting the y line text
axis y discontinuity=crunch, enlargelimits=false,
%axis x line*=bottom,
axis y line*=left,
ytick={-32,-28,-24, -20, -16},
ylabel style={font=\color{mycolor1}},
ylabel={$\Delta J_{d}$ [\%] },
axis background/.style={fill=white}
]
\addplot [color=mycolor1, line width=1.4pt, forget plot]
  table[row sep=crcr]{%
  	50	-18.9722809467604\\
  	80	-28.9944299948698\\
  	90	-30.6397649348534\\
  	100	-32.0676844026196\\
  	110	-32.9207146972765\\
  	120	-33.4509072736196\\
  	130	-33.4602089093844\\
  	140	-33.407761822074\\
  	150	-33.4755654948335\\
  	160	-33.2433200581226\\
  	170	-33.1361196422544\\
  	180	-33.106460095258\\
  	190	-33.0630597757406\\
  	200	-33.010554214874\\
  	210	-32.9521225730841\\
  	220	-32.865071128008\\
  	230	-32.8050275054911\\
  	240	-32.7562537987652\\
  	250	-32.7749050711953\\
  	300	-32.7244804994915\\
  };
\end{axis}

\begin{axis}[%
width=\figurewidth,
height=\figureheight,
at={(0\figurewidth,0\figureheight)},
scale only axis,
xmin=50,
xmax=300,
every outer y axis line/.append style={mycolor2},
every y tick label/.append style={font=\color{mycolor2}},
every y tick/.append style={mycolor2},
ymin=0,
ymax=2,
ytick={0, 0.5, 1, 1.5, 2},
xtick = {50, 100, 150, 200, 250, 300},
ylabel style={font=\color{mycolor2}},
xlabel={$T$ [s]},
ylabel={$\bar{t}_{\mathrm{imp}}$ [s]},
axis x line*=bottom,
axis y line*=right,
ylabel style={yshift=0.2cm}, %shifting the y line text
]
\addplot[area legend, draw=none, fill=white!80!red, fill opacity=0.8, forget plot]
table[row sep=crcr] {%
	x	y\\
	50	0.158631563186646\\
	80	0.238743233680725\\
	90	0.281968486309052\\
	100	0.295747256278992\\
	110	0.324932098388672\\
	120	0.35006650686264\\
	130	0.379831075668335\\
	140	0.403641092777252\\
	150	0.430657649040222\\
	160	0.456530523300171\\
	170	0.490134215354919\\
	180	0.516068196296692\\
	190	0.546606969833374\\
	200	0.580662405490875\\
	210	0.607179117202759\\
	220	0.650837790966034\\
	230	0.678898096084595\\
	240	0.711255037784576\\
	250	0.731899118423462\\
	300	1.10942158699036\\
	300	1.79820532798767\\
	250	1.95882420539856\\
	240	1.5942519903183\\
	230	1.55897221565247\\
	220	1.42815897464752\\
	210	1.41077835559845\\
	200	1.39606431722641\\
	190	1.35139517784119\\
	180	1.26843711137772\\
	170	1.2101354598999\\
	160	1.15551176071167\\
	150	1.13250865936279\\
	140	1.11912138462067\\
	130	1.10208640098572\\
	120	1.23574004173279\\
	110	1.04239416122437\\
	100	1.3038782119751\\
	90	0.966308581829071\\
	80	0.9852388381958\\
	50	0.853491067886353\\
	50	0.158631563186646\\
}--cycle;
\addplot [color=red, line width=1.4pt, forget plot]
  table[row sep=crcr]{%
  	50	0.179235100746155\\
  	80	0.264176249504089\\
  	90	0.310412645339966\\
  	100	0.318177223205566\\
  	110	0.368820548057556\\
  	120	0.3872230052948\\
  	130	0.413279891014099\\
  	140	0.447013139724731\\
  	150	0.542755126953125\\
  	160	0.561555504798889\\
  	170	0.613661050796509\\
  	180	0.676716327667236\\
  	190	0.718522071838379\\
  	200	0.784644365310669\\
  	210	0.812060236930847\\
  	220	0.881754636764526\\
  	230	0.909811735153198\\
  	240	0.968738317489624\\
  	250	0.990082740783691\\
  	300	1.17715096473694\\
  };
\end{axis}

\end{tikzpicture}% 
	\caption{ \label{fig:rhimp} \small  Results from the improvement step when planning a motion to $\xbm_{f,1}$ for different planning horizons $T$. Blue: improvement in objective function value compared to the nominal trajectory. Red: median computation time for the improvement step. The shaded area shows the span of computation times below the 95th~percentile.  }
\end{figure}

First, the impact of using different planning horizons in the improvement step is analyzed. The desired sampling time is selected as $\Delta_t = 2$ s in \eqref{eq:ocp2}. Hence, the number of discretization points $N$ defines the desired planning horizon $T=\Delta_t N$. A suitable value of $T$ is both system and problem specific. Here, the value of $T$ is selected by evaluating several planning horizons applied on the scenarios shown in Fig.~\ref{fig:sc}. The results are summarized in Fig.~\ref{fig:rhimp}, which presents the reduction in objective function value compared to the nominal trajectory as a function of the planning horizon $T$. Furthermore, the computation times for the corresponding values of $T$ are also displayed. The results indicate that after a planning horizon of 150 s, no further improvement in objective function value is obtained, while the computation time continues to increase. Therefore, $T =150$ s is a reasonable choice and is used in the remainder of this section. Moreover, with $T=150$ s, the 95th percentile of the computation time is 1.13 s (median: 0.7 s), which is well below the desired sampling time \mbox{$\Delta_t = 2$ s}. 

\begin{figure}
	\centering
	\setlength\figureheight{0.227\textwidth}
	\setlength\figurewidth{0.4\textwidth}
	\input{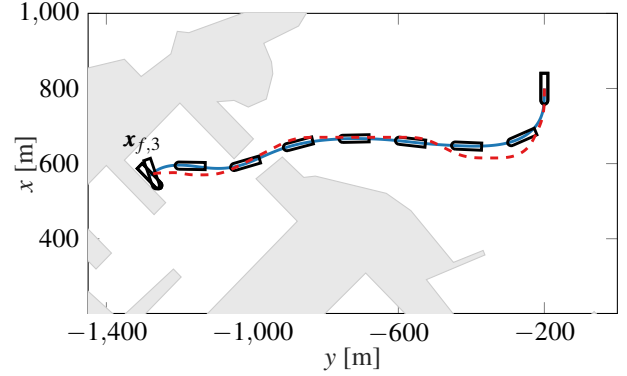} \vspace*{-0.5em}
	\caption{ \label{fig:ship_path} \small The nominal trajectory (dotted red) and the solution from the improvement step (blue) when planning a motion to $\xbm_{f,3}$. \vspace*{-0.5em} }
\end{figure}

\begin{figure}[]
	\centering \hspace*{0.5em} \vspace*{-0.5em}
	\subfloat[Body-fixed velocities (m/s)\label{subfig-1:dummy}]{%
		\setlength\figureheight{0.15\textwidth}
		\setlength\figurewidth{0.4\textwidth}
		\input{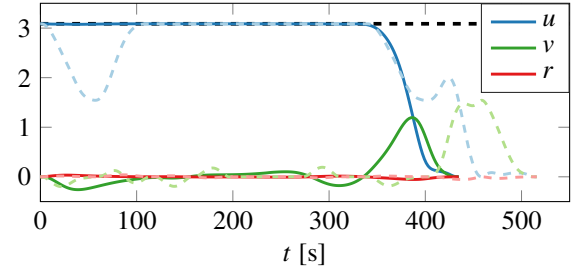}
	} \\ \vspace*{-0.5em}
	\subfloat[Thruster angles ($^\circ$) \label{subfig-2:dummy}]{%
		\setlength\figureheight{0.15\textwidth}
		\setlength\figurewidth{0.4\textwidth}
		\input{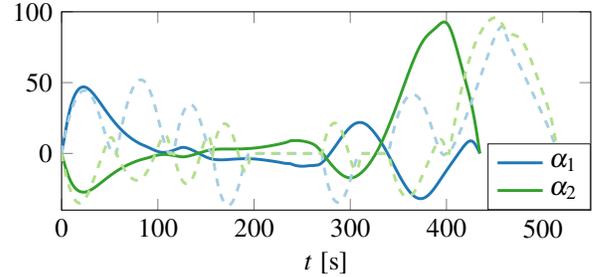}
	}    \\ \hspace*{0.5em} \vspace*{-0.5em}
	\subfloat[Thruster velocities (RPS)\label{subfig-3:dummy}]{%
		\setlength\figureheight{0.15\textwidth}
		\setlength\figurewidth{0.4\textwidth}
		\input{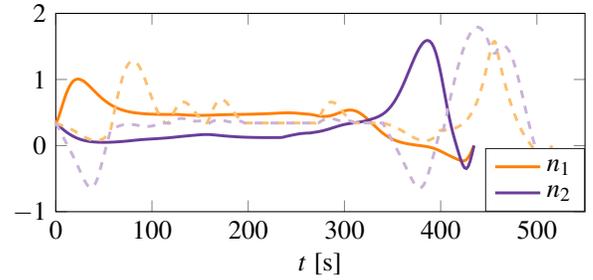}
	}
	\caption{ \small Trajectories for the problem illustrated in Fig.~\ref{fig:ship_path} using a planning horizon $T=150$ s. The dashed lines represent the nominal trajectory from the lattice-based planner. }
	\label{fig:dummy}
\end{figure}

Table~\ref{tab:res_summary} highlights the benefits of improving the nominal trajectory.
Both the time $T_f$ to reach the terminal state and
the amount of force applied by the thrusters $F_{\mathrm{tot}}$ are reduced.  
In particular, the applied force is significantly reduced by 40--50$\%$ in the given scenarios.  The enhancements come at the expense of a minor latency time, arising from the computation time of the initial spatial constraints $t_{\mathrm{env}}$ and the first receding-horizon improvement step. These minor initial computations could, however, be performed already when the ship is approaching the initial state at which the plan starts. In summary, smoother and more energy-efficient solutions are in these scenarios obtained using the receding-horizon improvement step. These results are illustrated for one scenario in Fig.~\ref{fig:ship_path}--\ref{fig:dummy}, where the trajectory computed by the improvement step is compared to the nominal trajectory computed by the lattice-based planner.

\begin{figure}
	\centering
	\setlength\figureheight{0.12\textwidth}
	\setlength\figurewidth{0.4\textwidth}
	\input{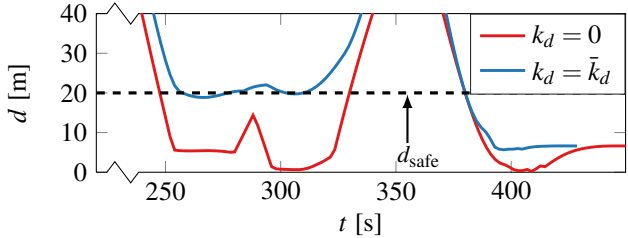}
	\caption{ \label{fig:dist} \small Distance to obstacles along the trajectory for the problem illustrated in Fig.~\ref{fig:ship_path}. The blue line is obtained when using a value of $k_d = \bar{k}_d = 1.5\cdot 10^{-3}$, and the red line when $k_d = 0$. The black dashed line represents the desired safety distance $d_{\mathrm{safe}}$. }
\end{figure}

Finally, Fig.~\ref{fig:dist} illustrates the impact of adding the safety-distance parameter to~\eqref{eq:ocp2}. When the distance to obstacles is not penalized, i.e., $k_d = 0$, it can be seen that the ship is very close to the boundary of an obstacle ($d_{\text{min}} = 0.2$ m). In contrast, when $k_d = \bar{k}_d = 1.5\cdot 10^{-3}$, the desired safety distance $d_{\mathrm{safe}}$ is only violated during the final phase of the maneuver \mbox{($t > 370$ s)} mainly due to the location of the terminal state. Note that the terminal state is reached \emph{faster} when $k_d > 0$. This is because the environmental cost is integrated over time, which means that the time spent with $d(t) < d_\mathrm{safe}$ is minimized. 
One possibility to alleviate this behavior is to choose the value of $k_d$ adaptively based on, e.g., the current distance to the terminal state.  

\begin{table*}[]
	\caption{\small Results from the four docking problems in Fig.~\ref{fig:sc} using a planning horizon $T = 150$ s. $t_{\mathrm{lat}}$ and $t_{\mathrm{env}}$ are the times to compute the nominal trajectory and the initial spatial envelopes, respectively. $\bar{t}_{\text{imp}}$ is the median computation time for the improvement step. Finally $T_f$ and $F_{\mathrm{tot}}$ are the terminal time and total applied force, respectively, for the nominal trajectory (lat) and after the improvement step (imp).    } \label{tab:res_summary}	
	\normalsize
	\centering 
	\setlength{\tabcolsep}{4pt}
	\begin{tabular}{cccccccc}	
		$\xbm_{f,j}$ & $t_{\mathrm{lat}}$ [s] & $t_{\mathrm{env}}$ [s] & $\bar{t}_{\text{imp}}$ [s] & $T_{f, \mathrm{lat}}$ [s] & $T_{f, \mathrm{imp}}$ [s] &$F_{\mathrm{tot}, \mathrm{lat}}$ [MN] & $F_{\mathrm{tot}, \mathrm{imp}}$ [MN]   \\
		\hline
		1  & 3.0 & 1.4 & 0.46 & 652 & 617 & $152$  & 90   \\
		2  & 0.19 & 1.3 & 0.55 & 563 & 532 & 131 & 67 \\
		3 & 2.3 & 1.1 & 0.54 & 510 & 435  & 155 & 86  \\ 
		4 & 6.6 & 1.1 & 0.44 & 826 & 807  & 206 & 120  \\ 
		\hline
	\end{tabular}
	\vspace*{-1.5em}
\end{table*}
\vspace*{-0.5em}
\section{Conclusions and Future Work}\label{sec:conclusions}
This paper introduces a two-step optimization-based motion planner for maneuvering of ships in confined environments such as harbors. In the first step, a lattice-based motion planner solves the combinatorial aspects of the problem and computes a feasible, but suboptimal solution to the problem using a finite library of optimized maneuvers. In the second step, the solution from the lattice-based planner is improved using receding-horizon optimization techniques. For efficient obstacle avoidance, a novel algorithm for automatic computation of safe spatial envelopes is introduced. Furthermore, a safety-distance parameter is introduced in the optimization-based improvement step which enables the algorithm to maintain a desired safety distance to obstacles without reducing the feasible region. The proposed motion planner is applied to a set of autonomous docking problems in a model of the Cape Town harbor, where it successfully computes energy-efficient solutions even in narrow passages.

Future work includes to extend the algorithm such that it is possible to handle dynamic scenarios including other ships. Another interesting topic for future work is to investigate how to incorporate regulations such as COLREGs~\cite{kufoalor2019autonomous} in the motion planning framework.
\vspace*{-0.8em}
\newpage
\bibliography{myrefs.bib}		

\begin{thebibliography}{29}
\providecommand{\natexlab}[1]{#1}
\providecommand{\url}[1]{\texttt{#1}}
\expandafter\ifx\csname urlstyle\endcsname\relax
  \providecommand{\doi}[1]{doi: #1}\else
  \providecommand{\doi}{doi: \begingroup \urlstyle{rm}\Url}\fi

\bibitem[Fossen(2011)]{fossen2011handbook}
Thor~I Fossen.
\newblock \emph{Handbook of marine craft hydrodynamics and motion control}.
\newblock John Wiley \& Sons, 2011.

\bibitem[Perez(2006)]{perez2006ship}
Tristan Perez.
\newblock \emph{Ship motion control: course keeping and roll stabilisation
  using rudder and fins}.
\newblock Springer Science \& Business Media, 2006.

\bibitem[Veksler et~al.(2016)]{veksler2016dynamic}
Aleksander Veksler et~al.
\newblock Dynamic positioning with model predictive control.
\newblock \emph{IEEE Trans. Control Syst. Technol.}, 24\penalty0 (4):\penalty0
  1340--1353, 2016.

\bibitem[Moreira et~al.(2007)Moreira, Fossen, and Soares]{moreira2007path}
L{\'u}cia Moreira, Thor~I Fossen, and C~Guedes Soares.
\newblock Path following control system for a tanker ship model.
\newblock \emph{Ocean Engineering}, 34\penalty0 (14-15):\penalty0 2074--2085,
  2007.

\bibitem[Martinsen et~al.(2019)Martinsen, Lekkas, and
  Gros]{martinsen2019autonomous}
Andreas~B Martinsen, Anastasios~M Lekkas, and Sebastien Gros.
\newblock Autonomous docking using direct optimal control.
\newblock In \emph{Proc. 12th IFAC Conf. Control Appl. in Marine Systems,
  Robotics, and Veh. (CAMS)}, volume~52, pages 97--102. IFAC, 2019.

\bibitem[Du et~al.(2018)]{du2018motion}
Zhe Du et~al.
\newblock Motion planning for unmanned surface vehicle based on trajectory
  unit.
\newblock \emph{Ocean Engineering}, 151:\penalty0 46--56, 2018.

\bibitem[Han et~al.(2020)]{han2020autonomous}
Jungwook Han et~al.
\newblock Autonomous collision detection and avoidance for {ARAGON USV}:
  Development and field tests.
\newblock \emph{J. Field. Robot.}, 2020.

\bibitem[Kufoalor et~al.(2019)]{kufoalor2019autonomous}
DKM Kufoalor et~al.
\newblock Autonomous maritime collision avoidance: Field verification of
  autonomous surface vehicle behavior in challenging scenarios.
\newblock \emph{J. Field. Robot.}, 2019.

\bibitem[Bergman et~al.(June 2019)Bergman, Ljungqvist, and
  Axehill]{bergman2019improved}
Kristoffer Bergman, Oskar Ljungqvist, and Daniel Axehill.
\newblock Improved optimization of motion primitives for motion planning in
  state lattices.
\newblock In \emph{2019 IEEE Intell. Veh. Symp. (IV)}, June 2019.

\bibitem[Bergman et~al.({\natexlab{a}})Bergman, Ljungqvist, and
  Axehill]{bergman2019bimproved}
Kristoffer Bergman, Oskar Ljungqvist, and Daniel Axehill.
\newblock Improved path planning by tightly combining lattice-based path
  planning and numerical optimal control.
\newblock \emph{Accepted for publ. in IEEE Trans. Intell. Veh.}. Pre-print
  available at {arXiv}: \url{https://arxiv.org/abs/1903.07900}, {\natexlab{a}}.

\bibitem[Bergman et~al.({\natexlab{b}})]{bergman2020optimization}
Kristoffer Bergman et~al.
\newblock An optimization-based receding horizon planning algorithm.
\newblock \emph{Accepted for publ. at 21st IFAC World Congress}, 2020,
  {\natexlab{b}}.

\bibitem[Pivtoraiko et~al.(2009)Pivtoraiko, Knepper, and
  Kelly]{pivtoraiko2009differentially}
Mihail Pivtoraiko, Ross~A Knepper, and Alonzo Kelly.
\newblock Differentially constrained mobile robot motion planning in state
  lattices.
\newblock \emph{J. Field. Robot.}, 26\penalty0 (3):\penalty0 308--333, 2009.

\bibitem[Andersson et~al.(2018)]{LjungqvistCDC2018}
Olov Andersson et~al.
\newblock Receding-horizon lattice-based motion planning with dynamic obstacle
  avoidance.
\newblock In \emph{Proc. 57th IEEE Conf. Decis. Control}, 2018.

\bibitem[Ljungqvist et~al.(2019)]{LjungqvistJFR2019}
O.~Ljungqvist et~al.
\newblock A path planning and path-following control framework for a general
  2-trailer with a car-like tractor.
\newblock \emph{J. Field. Robot.}, 36\penalty0 (8):\penalty0 1345--1377, 2019.

\bibitem[Likhachev and Ferguson(2009)]{likhachev2009planning}
M.~Likhachev and D.~Ferguson.
\newblock Planning long dynamically feasible maneuvers for autonomous vehicles.
\newblock \emph{The Int. J. Robot. Res.}, 28\penalty0 (8):\penalty0 933--945,
  2009.

\bibitem[Cirillo et~al.(2014)Cirillo, Uras, and Koenig]{cirillo2014lattice}
Marcello Cirillo, Tansel Uras, and Sven Koenig.
\newblock A lattice-based approach to multi-robot motion planning for
  non-holonomic vehicles.
\newblock In \emph{{Internat.} Conf. on Intell. Robots and Systems (IROS)},
  pages 232--239. IEEE, 2014.

\bibitem[Andreasson et~al.(2015)]{andreasson2015fast}
Henrik Andreasson et~al.
\newblock Fast, continuous state path smoothing to improve navigation accuracy.
\newblock In \emph{IEEE Internat. Conf. on Robot. and Autom. (ICRA)}, pages
  662--669. IEEE, 2015.

\bibitem[Liu et~al.(2017)]{liu2017planning}
Sikang Liu et~al.
\newblock Planning dynamically feasible trajectories for quadrotors using safe
  flight corridors in 3-d complex environments.
\newblock \emph{IEEE Robot. and Autom. Lett.}, 2\penalty0 (3):\penalty0
  1688--1695, 2017.

\bibitem[Schoels et~al.()]{schoels2020ciao}
Tobias Schoels et~al.
\newblock {CIAO}$^\star$: {MPC}-based safe motion planning in predictable
  dynamic environments.
\newblock \emph{Accepted for publ. at 21st IFAC World Congress}, 2020.

\bibitem[LaValle(2006)]{lavalle2006planning}
Steven~M LaValle.
\newblock \emph{Planning Algorithms}.
\newblock Cambridge University Press, Cambridge, UK, 2006.

\bibitem[Lewandowski(2004)]{Lewandowski}
Edward~M Lewandowski.
\newblock \emph{The dynamics of marine craft: maneuvering and seakeeping},
  volume~22.
\newblock World scientific, 2004.

\bibitem[Pivtoraiko and Kelly(2005)]{pivtoraiko2005efficient}
Mihail Pivtoraiko and Alonzo Kelly.
\newblock Efficient constrained path planning via search in state lattices.
\newblock In \emph{Internat. Sym. on Artificial Intell., Robot., and Autom. in
  Space}, pages 1--7, 2005.

\bibitem[Ferguson and Likachev(2008)]{ferguson2008efficient}
Dave Ferguson and Maxim Likachev.
\newblock Efficiently using cost maps for planning complex maneuvers.
\newblock In \emph{Proc. of Internat. Conf. on Robot. and Autom. Workshop on
  Planning with Cost Maps}. IEEE, 2008.

\bibitem[Knepper and Kelly(2006)]{knepper2006high}
Ross~A Knepper and Alonzo Kelly.
\newblock High performance state lattice planning using heuristic look-up
  tables.
\newblock In \emph{2006 IEEE/RSJ Inter. Conf. on Intell. Robots and Systems},
  pages 3375--3380, 2006.

\bibitem[Diehl et~al.(2006)]{diehl2006fast}
Moritz Diehl et~al.
\newblock Fast direct multiple shooting algorithms for optimal robot control.
\newblock In \emph{Fast motions in biomechanics and robotics}, pages 65--93.
  Springer, 2006.

\bibitem[B{\"a}rlund et~al.()]{barlund2020nonlinear}
Alexander B{\"a}rlund et~al.
\newblock Nonlinear {MPC} for combined motion control and thrust allocation of
  ships.
\newblock \emph{Accepted for publ. at 21st IFAC World Congress}, 2020.

\bibitem[Andersson et~al.(In Press, 2018)]{andersson2018casadi}
Joel A~E Andersson et~al.
\newblock {CasADi} -- {A} software framework for nonlinear optimization and
  optimal control.
\newblock \emph{Mathematical Programming Computation}, In Press, 2018.

\bibitem[W{\"a}chter and Biegler(2006)]{wachter2006implementation}
Andreas W{\"a}chter and Lorenz~T Biegler.
\newblock On the implementation of an interior-point filter line-search
  algorithm for large-scale nonlinear programming.
\newblock \emph{Mathematical programming}, 106\penalty0 (1):\penalty0 25--57,
  2006.

\bibitem[Perez et~al.(2006)]{perez2006overview}
Tristan Perez et~al.
\newblock An overview of the marine systems simulator ({MSS}): A simulink
  toolbox for marine control systems.
\newblock \emph{Modeling, Identification and Control}, 27\penalty0
  (4):\penalty0 259--275, 2006.

\end{thebibliography}
\end{document}